\documentclass[letterpaper, 10pt, conference]{IEEEconf}
\IEEEoverridecommandlockouts 

\usepackage{amsmath,amssymb}
\usepackage{graphicx,subfigure,url,array,threeparttable}

\newcommand{\norm}[1]{\ensuremath{\left\| #1 \right\|}}
\newcommand{\bracket}[1]{\ensuremath{\left[ #1 \right]}}
\newcommand{\braces}[1]{\ensuremath{\left\{ #1 \right\}}}
\newcommand{\parenth}[1]{\ensuremath{\left( #1 \right)}}
\newcommand{\refeqn}[1]{(\ref{eqn:#1})}
\newcommand{\reffig}[1]{Fig. \ref{fig:#1}}
\newcommand{\tr}[1]{\mbox{tr}\ensuremath{\negthickspace\bracket{#1}}}

\newcommand{\SO}{\ensuremath{\mathrm{SO(3)}}}
\newcommand{\T}{\ensuremath{\mathrm{T}}}
\newcommand{\so}{\ensuremath{\mathfrak{so}(3)}}

\renewcommand{\Re}{\ensuremath{\mathbb{R}}}
\renewcommand{\S}{\ensuremath{\mathbb{S}}}

\title{\LARGE \bf
Optimal Attitude Control for a Rigid Body with Symmetry}

\author{ \parbox{3 in}{\centering Taeyoung Lee\authorrefmark{1}\authorrefmark{2}, N. Harris McClamroch\authorrefmark{2}\\
         Department of Aerospace Engineering\\
         University of Michigan, Ann Arbor, MI 48109\\
         {\tt\small \{tylee, nhm\}@umich.edu}}
         \hspace*{ 0.5 in}
         \parbox{3 in}{\centering Melvin Leok\authorrefmark{1}\\
         Department of Mathematics\\
        Purdue University, West Lafayette, IN 47907\\
         {\tt\small mleok@math.purdue.edu}}
        \thanks{\textsuperscript{\footnotesize\ensuremath{*}}This research has been supported in part by NSF under grant DMS-0504747, and by a grant from the Rackham Graduate School, University of Michigan.}
        \thanks{\textsuperscript{\footnotesize\ensuremath{\dagger}}This research has been supported in part by NSF under grant ECS-0244977.}
}

\begin{document}
\allowdisplaybreaks
\maketitle \thispagestyle{empty} \pagestyle{empty}

\begin{abstract}
Optimal control problems are formulated and efficient computational
procedures are proposed for attitude dynamics of a rigid body with
symmetry. The rigid body is assumed to act under a gravitational
potential and under a structured control moment that respects the
symmetry. The symmetry in the attitude dynamics system yields a
conserved quantity, and it causes a fundamental singularity in the
optimal control problem. The key feature of this paper is its use of
computational procedures that are guaranteed to avoid the numerical
ill-conditioning that originates from this symmetry. It also
preserves the geometry of the attitude dynamics. The theoretical
basis for the computational procedures is summarized, and examples
of optimal attitude maneuvers for a 3D pendulum are presented.
\end{abstract}

\section{Introduction}

We study a discrete optimal control problem for attitude dynamics of
a rigid body with symmetry. The attitude is represented by a
rotation matrix, which has a Lie group structure denoted by \SO. We
assume that the rigid body is acting under an attitude dependent
potential, and the potential is invariant under a symmetry action.
The external control input is formulated such that it respects the
symmetry.  This problem provides both a theoretical challenge and a
numerical challenge in the sense that the configuration space has a
Lie group structure, and the conserved quantity causes
ill-conditioning of the numerical optimization.

General purpose numerical integration methods, including the popular
Runge-Kutta schemes, typically preserve neither the group structure
of the attitude configuration space nor the invariant properties of
the dynamics. Geometric structure-preserving integrators are
symplectic and momentum preserving, and they exhibit good energy
behavior for an exponentially long time
period~\cite{HaiLubWan.BK02}. In particular, Lie group variational
integrators have the desirable properties that they preserve the
group structure as well as the geometric features, without needs of
local parameterization, reprojection, or
constraints~\cite{CCA05,CMAME05}. The exact geometric properties of
the discrete flow not only generate improved qualitative behavior,
but also allow for accurate long-time simulation.

Optimal control problems on a Lie group have been studied in
\cite{Spin.MCSS98,Sas.ICIAM95,Jur.BK97}. These studies are based on
the driftless kinematics of a Lie group. The dynamics are ignored,
and elements in the corresponding Lie algebra are considered as
control inputs. The discrete optimal control problems of a rigid
body are studied in \cite{JOTA06,CDC06.1}, where the dynamics as
well as the kinematics equations are explicitly utilized, and an
efficient numerical algorithm to solve discrete optimality
conditions is presented.

This paper introduces geometrically exact and numerically efficient
computational approaches to solve the optimal control problems of
the attitude dynamics of a rigid body with symmetry and structured
control input. The dynamics are discretized by a Lie group
variational integrator, and discrete necessary conditions for
optimality are constructed. The utilization of the Lie group
variational integrator is justified in this problem, since it
preserves the momentum map originating from the symmetry. The rigid
body is underactuated since the control input does not act along the
symmetry direction. The symmetry of the controlled dynamics causes
difficulties in solving the necessary conditions for optimality. A
simple numerical approach is presented to overcome this numerical
ill-conditioning.

This paper is organized as follows. In Section \ref{sec:dyn}, a 3D
pendulum is presented as a model of rigid body attitude dynamics,
and the symmetry of the 3D pendulum is described. An optimal control
problem with symmetry is studied in Section \ref{sec:opt}, and
numerical results are given in Section \ref{sec:ne}.

\section{Dynamics of a 3D Pendulum}\label{sec:dyn}
A 3D pendulum is a rigid body supported by a fixed frictionless
pivot acting under the influence of uniform gravitational
field~\cite{SheSanNalBerMcC.CDC04}. We use a 3D pendulum model to
study the optimal control for attitude dynamics of a rigid body,
since it has three degrees of rotational freedom, and the
gravitational potential has a symmetry: it is invariant under a
rotation about the gravity direction.

In this section, the continuous equations of motion are presented.
The symmetry of the 3D pendulum are discussed, and the control input
structure is described. Discrete equations of motion, referred to as
a Lie group variational integrator, are described for a controlled
3D pendulum model.

\subsection{Continuous equations of motion}
The configuration space of the 3D pendulum is $\SO$. We identify the
tangent bundle $\T\SO$ with $\SO\times\so$ by left translation, and
we identify $\so$ with $\Re^3$ by an isomorphism
$S(\cdot):\Re^3\mapsto\so$ defined by the condition that $S(x)y=x
\times y$ for any $x,y\in\Re^3$. We denote the attitude and the
angular velocity of the rigid body as $(R,\Omega)\in\T_R\SO$. The
rotation matrix $R\in\SO$ transforms a vector represented in the
body fixed frame to one represented in the inertial frame.

Let $\rho\in\Re^3$ be a vector from the pivot point to the mass
center of the rigid body expressed in the body fixed frame, and let
$m,g\in\Re$ and $J\in\Re^{3\times 3}$ be the mass of the rigid body,
the gravitational acceleration, and the moment of inertia matrix of
the rigid body about the pivot point, respectively. The Lagrangian
of the 3D pendulum $L:\T\SO\mapsto\Re$ is given by
\begin{align*}
    L(R,\Omega)=\frac{1}{2}\tr{S(\Omega) J_d S(\Omega)^T} + mg e_3^T R \rho,
\end{align*}
where $J_d\in\Re^{3\times 3}$ is a nonstandard moment of inertia
defined by $J_d=\frac{1}{2}\tr{J}I_{3\times 3}-J$, and we set the
gravitational direction in the inertial frame as
$e_3=[0;0;1]\in\S^2$.

The continuous equations, derived from the Lagrange-d'Alembert
principle, are given by
\begin{gather}
    \dot \Pi + \Omega\times\Pi = mg \rho \times R^T e_3+M,\label{eqn:Pidot}\\
    \dot R = R S(\Omega),\label{eqn:Rdot}
\end{gather}
where $\Pi=J\Omega\in\Re^3$ is the angular momentum in the body
fixed frame, and $M\in\Re^3$ is the external control moment.

\subsection{Symmetry of 3D pendulum}
The kinetic energy of the rigid body is left invariant on $\T\SO$,
and the gravitational potential energy is invariant under a rotation
about the gravity direction, which can be represented by the left
action of the subgroup $\braces{\exp S(\theta e_3)\in\SO\big|
\theta\in\S^1}$.

As a result, the Lagrangian of the 3D pendulum has a symmetry action
by $\S^1$, $\Phi_\theta:\S^1\times \SO\mapsto\SO$ given by
\begin{align*}
    \Phi_\theta(R)=\exp S(\theta e_3)\, R,
\end{align*}
for $\theta\in\S^1$. It can be shown that $\Phi_\theta^*
L(R,\Omega)=L(R,\Omega)$.

Suppose that there is no external control input. Noether's theorem
states that a symmetry in the Lagrangian yields conservation of the
momentum map. For the 3D pendulum, the momentum map of the symmetry
action $\Phi_\theta$ corresponds to the inertial angular momentum of
the rigid body about the gravity direction $\pi_3=e_3^T
RJ\Omega\in\Re$. It is conserved for the free dynamics of the 3D
pendulum.

The structure of the control input respects the symmetry of the
uncontrolled free dynamics of the 3D pendulum, namely
\begin{align*}
    M=R^T e_3 \times u,
\end{align*}
for a control parameter $u\in\Re^3$. Since the external control
moment has no component along the gravity direction, the angular
momentum about the gravity direction is also preserved in the
controlled dynamics. Such control inputs are physically utilized by
actuation mechanisms, such as point mass actuators, that change the
center of mass of the 3D pendulum.

Here we introduce the concept of a geometric phase, and it is used
to interpret the numerical optimization result in Section
\ref{sec:ne}. Using the symmetry, the dynamics of the 3D pendulum
can be expressed in terms of $\Gamma=R^T e_3$ in the reduced
configuration space $\SO/\S^1\simeq \S^2$. The corresponding flow in
the original configuration space is reconstructed by lifting to a
level set of the conserved quantity. Suppose that the trajectory in
the reduced space is a closed loop, i.e. $\Gamma(0)=\Gamma(T)$ for
some $T>0$, and the value of the angular momentum about the gravity
direction is zero. Then, the terminal attitude is related to the
initial attitude by a symmetric action. More explicitly, we have
\begin{align*}
    R(T)=\Phi_{\theta_{\text{geo}}} (R(0)),
\end{align*}
where $\theta_{\text{geo}}$ is the geometric phase determined by
\begin{align}
    \theta_{\text{geo}}=\int_{\mathcal{B}}
    \frac{2\norm{J\Gamma}^2-\tr{J}(\Gamma^T J\Gamma)}{(\Gamma^T
    J\Gamma)^2}\,dA,\label{eqn:geo}
\end{align}
where $\mathcal{B}$ is a surface in $\S^2$ whose boundary is
$\braces{\Gamma(t)\big| t\in[0,T]}$~\cite{MarMonRat.BK90}. Note that
the geometric phase is determined only by the reduced trajectory of
$\Gamma$ and the characteristics of the rigid body $J$. It is
independent of the velocity $\dot\Gamma$.

\subsection{Lie group variational integrator}
The attitude of the 3D pendulum is represented by a rotation matrix
$R\in\SO$. The conserved quantity, arising from symmetry, is
emphasized in this study. However, the most common numerical
integration methods, including the widely used Runge-Kutta schemes,
neither preserve the Lie group structure nor first integrals. In
addition, standard Runge-Kutta methods fail to capture the energy
dissipation of a controlled system accurately~\cite{MarWes.AN01}.
For example, if we integrate \refeqn{Rdot} by a typical Runge-Kutta
scheme, the quantity $R^T R$ inevitably drifts from the identity
matrix as the simulation time increases. It is often proposed to
parameterize \refeqn{Rdot} by Euler angles or unit quaternions.
However, Euler angles are not global expressions of the attitude
since they have associated singularities. Unit quaternions do not
exhibit singularities, but are constrained to lie on the unit
three-sphere $\S^3$, and general numerical integration methods do
not preserve the unit length constraint. Therefore, quaternions have
the same numerical drift problem. Renormalizing the quaternion
vector at each step tends to break other conservation properties.
Furthermore, unit quaternions, which are diffeomorphic to
$\mathrm{SU(2)}$, double cover $\SO$. So there are inevitable
ambiguities in expressing the attitude using quaternions.

In~\cite{CCA05}, Lie group variational integrators are introduced by
explicitly adopting the approach of Lie group
methods~\cite{IseMunHan.AN00} to the discrete variational
principle~\cite{MarWes.AN01}. They have the desirable property that
they are symplectic and momentum preserving, and they exhibit good
energy behavior for an exponentially long time period. They also
preserve the Lie group structure without the use of local charts,
reprojection, or constraints.

Using the results in~\cite{CCA05}, a Lie group variational
integrator on $\SO$ is given for the 3D pendulum by
\begin{gather}
h S(\Pi_k) = F_k J_d - J_dF_k^T,\label{eqn:findf}\\
R_{k+1} = R_k F_k,\label{eqn:updateR}\\
\Pi_{k+1} = F_k^T \Pi_k + hmg\rho\times R_{k+1}^T e_3+hR_{k+1}^T e_3
\times u_{k+1},\label{eqn:updatePi}
\end{gather}
where the subscript $k$ denotes the $k$th discrete variable for a
fixed integration step size $h\in\Re$, and $F_k\in\SO$ is the
relative attitude between two adjacent integration steps. For a
given $(R_k,\Pi_k)$ and control inputs, \refeqn{findf} is solved to
find $F_k$. Then $(R_{k+1},\Pi_{k+1})$ is obtained by
\refeqn{updateR} and \refeqn{updatePi}. This yields a map
$(R_k,\Pi_k)\mapsto(R_{k+1},\Pi_{k+1})$ and this process is
repeated. The only implicit part is \refeqn{findf}. The actual
computation of $F_k$ is done in the Lie algebra $\so$ of dimension
3, and the rotation matrices are updated by multiplication. So this
approach is distinguished from integration of the kinematics
equation \refeqn{Rdot}, and there is no excessive computational
burden. The properties of these discrete equations of motion are
discussed more explicitly in \cite{CCA05,CMAME05}. We use these
discrete equations of motion to formulate the following optimal
control problem.

\section{Optimal control with symmetry}\label{sec:opt}
We formulate an optimal attitude control problem for a 3D pendulum
with symmetry. Necessary conditions for optimality are developed and
computational approaches are presented to solve the corresponding
two point boundary value problem.

\subsection{Problem formulation}

A discrete time optimal control problem is formulated as a maneuver
of the rigid pendulum body from a given initial attitude $R_0\in\SO$
and an initial angular momentum $\Pi_0\in\Re^3$ to a desired
terminal attitude $R_N^d\in\SO$ and a terminal angular momentum
$\Pi_N^d\in\Re^3$ during a given maneuver time $N$. The performance
index is the square of the $l_2$ norm of the control inputs:
\begin{gather*}
\text{given: } (R_0,\Pi_0),\,(R_N^d,\Pi_N^d),\,N,\\
\min_{u_{k+1}} \mathcal{J}=\sum_{k=0}^{N-1} \frac{h}{2}\norm{u_{k+1}}^2,\\
\text{such that } R_N=R_N^d,\,\Pi_N=\Pi_N^d,\\
\text{subject to \refeqn{findf}, \refeqn{updateR} and
\refeqn{updatePi}.}
\end{gather*}

\subsection{Necessary conditions of optimality}

\paragraph*{Variational models} The necessary conditions of optimality are developed using the
standard variational approach. We first derive certain variational
formulas. The variation of $R_k\in\SO$ can be expressed in terms of
a Lie algebra element $S(\zeta_k)\in\so$ for $\zeta_k\in\Re^3$ and
the exponential map as
\begin{align*}
    R_k^{\epsilon} = R_k \exp \epsilon S(\zeta_k).
\end{align*}
The corresponding infinitesimal variation is given by
\begin{align}
    \delta R_k =\frac{d}{d\epsilon}\bigg|_{\epsilon=0}
    R_k \exp \epsilon S(\zeta_k)=R_k S(\zeta_k).\label{eqn:delRk}
\end{align}
This gives an expression for the infinitesimal variation of a Lie
group element in terms of its Lie algebra. Then, small perturbations
from a given trajectory can be written as
\begin{align}
    \Pi_k^\epsilon & = \Pi_k +\epsilon\delta\Pi_k,\label{eqn:Pike}\\
    R_k^\epsilon & = R_k + \epsilon R_k S(\zeta_k)
    +\mathcal{O}(\epsilon^2),\label{eqn:Rke}
\end{align}
where $\delta\Pi_k,\zeta_k$ are considered as elements of $\Re^3$.

We derive expressions for the constrained variation of $F_k$ using
its definition \refeqn{updateR} and the variation of the rotation
matrix \refeqn{Rke}. Since $F_k=R_k^T R_{k+1}$, the infinitesimal
variation $\delta F_k$ is given by
\begin{align*}
    \delta F_k & = \delta R_k^T R_{k+1}+R_k^T \delta R_{k+1},\\
    & = -S(\zeta_k)F_k +F_kS(\zeta_{k+1}).
\end{align*}
We can also write $\delta F_k = F_k S(\xi_k)$ for $\xi_k\in\Re^3$
using \refeqn{delRk}. Using the property $S(R^T x)=R^T S(x)R$ for
any $R\in\SO$ and $x\in\Re^3$, we obtain the constrained variation
of $F_k$ as
\begin{align}
    \xi_k = -F_k^T \zeta_k +\zeta_{k+1}.\label{eqn:xik}
\end{align}

We now relate the constrained variation of $\delta\Pi_k$ to $\xi_k$
by starting with \refeqn{findf}. Taking a variation of
\refeqn{findf}, we obtain
\begin{align*}
    h S(\delta\Pi_k) = F_k S(\xi_k) J_d + J_d S(\xi_k)F_k^T.
\end{align*}
Using the properties, $S(Rx)=RS(x)R^T$ and
$S(x)A+A^TS(x)=S(\braces{\tr{A}I_{3\times 3}-A}x)$ for any
$x\in\Re^3$, $A\in\Re^{3\times 3}$, and $R\in\SO$, the above
equation is rewritten as
\begin{align*}
    hS(\delta\Pi_k)=S(\braces{\tr{F_kJ_d}I_{3\times
    3}-F_kJ_d}F_k\xi_k).
\end{align*}
Thus, $\xi_k$ is given by
\begin{align}
    \xi_k=hF_k^T\braces{\tr{F_kJ_d}I_{3\times
    3}-F_kJ_d}^{-1}=\mathcal{B}_k\delta\Pi_k,\label{eqn:BBk}
\end{align}
where $\mathcal{B}\in\Re^{3\times 3}$.

\paragraph*{Necessary conditions}
Define an augmented performance index as
\begin{align}
\mathcal{J}_a = \sum_{k=0}^{N-1} & \frac{h}{2}\norm{u_{k+1}}^2
+\lambda_k^{1,T}S^{-1}\!\parenth{\mathrm{logm}(F_k-R_{k}^TR_{k+1})}\nonumber\\
& +\lambda_k^{2,T}\braces{-\Pi_{k+1} + F_k^T \Pi_k + hmg\rho\times
R_{k+1}^T e_3}\nonumber\\
& +\lambda_k^{2,T}\braces{R_{k+1}^T e_3 \times
u_{k+1}},\label{eqn:Ja}
\end{align}
where $\lambda_k^{1}, \lambda_k^{2}\in \Re^3$, are Lagrange
multipliers corresponding to the discrete equations of motion
\refeqn{updateR} and \refeqn{updatePi}. The constraint
\refeqn{findf} is applied implicitly by \refeqn{BBk} when taking the
variation.

Using the variational models \refeqn{Pike}--\refeqn{BBk}, and the
fact that the variations $\zeta_k,\,\delta\Pi_k$ vanish at
$k={0,N}$, the infinitesimal variation of the augmented performance
index is written as
\begin{align*}
\delta\mathcal{J}_a &= \sum_{k=1}^{N-1} h\delta u_{k}^T
\braces{u_{k}-R_{k}^T e_3\times \lambda^2_{k-1}}\nonumber\\
&+\zeta_k^T \braces{-\lambda_{k-1}^1+\mathcal{A}^T_k\lambda_k^1
+\mathcal{C}^T_k\lambda_k^2-hF_ku_{k+1}e_3^TR_{k+1}}\nonumber\\
&+\delta\Pi_k^T \braces{-\lambda_{k-1}^2+\mathcal{B}_k^T\lambda_k^1
+\mathcal{D}_k^T\lambda_k^2-h\mathcal{B}_k^Tu_{k+1}e_3^TR_{k+1}},
\end{align*}
where
\begin{align*}
\mathcal{A}_k & = F_k^T,\\
\mathcal{B}_k & = hF_k^T\braces{\tr{F_kJ_d}I_{3\times 3}-F_kJ_d}^{-1},\\
\mathcal{C}_k & = hmgS(\rho)S(R_{k+1}^Te_3)F_k^T,\\
\mathcal{D}_k & = F_k^T + S(F_k^T\Pi_k)\mathcal{B}_k+h
mgS(\rho)S(R_{k+1}^Te_3)\mathcal{B}_k.
\end{align*}

Since $\delta\mathcal{J}_a=0$ for all variations of $\delta u_k,
\zeta_k, \delta\Pi_k$, the expressions in the braces of the above
equation are zero. Thus we obtain necessary conditions for
optimality as follows.
\begin{gather}
h S(\Pi_k) = F_k J_d - J_dF_k^T,\label{eqn:findf1}\\
R_{k+1} = R_k F_k,\label{eqn:updateR1}\\
\Pi_{k+1} = F_k^T \Pi_k + hmg\rho\times R_{k+1}^T e_3+hR_{k+1}^T e_3
\times u_{k+1},\label{eqn:updatePi1}\\
u_{k+1} = R_{k+1}^T e_3\times \lambda^2_{k},\label{eqn:ukp}\\
\begin{bmatrix} \lambda_{k}^1 \\ \lambda_{k}^2 \end{bmatrix}
= \begin{bmatrix} \mathcal{A}_{k+1}^T & \mathcal{C}_{k+1}^T -hF_{k+1}u_{k+2}e_3^T R_{k+2}\\
\mathcal{B}_{k+1}^T & \mathcal{D}_{k+1}^T-h\mathcal{B}_{k+1}^Tu_{k+2}e_3^TR_{k+2}
\end{bmatrix}
\begin{bmatrix} \lambda_{k+1}^1 \\ \lambda_{k+1}^2\end{bmatrix}.\label{eqn:updatelam}
\end{gather}
In the above equations, the implicit parts are \refeqn{findf1} and
\refeqn{updatelam}. For a given initial condition
$(R_0,\Pi_0,\lambda_0^1,\lambda_0^2)$, we can find $F_0$ by solving
\refeqn{findf1}. Then, $R_1$ is obtained by \refeqn{updateR1}. Since
$u_1=R_1^Te_3\times \lambda_0^2$ by \refeqn{ukp}, $\Pi_1$ can be
obtained using \refeqn{updatePi1}. We solve \refeqn{findf1} to
obtain $F_1$ using $\Pi_1$. Finally, $\lambda_1^1,\lambda_1^2$ are
obtained by solving the implicit equation \refeqn{updatelam}, since
$\mathcal{A}_1,\mathcal{B}_1,\mathcal{C}_1,\mathcal{D}_1$ are
functions of $R_1,\Pi_1,F_1$.

The implicit equation \refeqn{findf1} is solved by Netwon's
iteration in the Lie algebra, and the implicit equation
\refeqn{updatelam} is solved by fixed point iteration. Numerical
computations show that two or three iterations are typically
required to achieve machine precision.

\subsection{Two point boundary value problem}\label{subsec:tpbvp}
The necessary conditions for optimality are given by a 12
dimensional two point boundary value problem. This problem is to
find the optimal discrete flow, multipliers, and control inputs to
satisfy the equations of motion \refeqn{findf1}--\refeqn{updatePi1},
optimality condition \refeqn{ukp}, multiplier equations
\refeqn{updatelam}, and boundary conditions simultaneously.

We substitute the optimality condition \refeqn{ukp} into the
equations of motion and the multiplier equations, and we apply the
shooting method to solve the two point boundary value problem using
sensitivity derivatives. The shooting method is numerically
efficient in the sense that the number of iteration parameters is
minimized; 6 elements of the initial Lagrange multiplier are
iterated. In other approaches, the entire discrete trajectory of the
control input and Lagrange multiplier are updated.

The drawback of the shooting method is that the extremal solutions
are sensitive to small changes in the unspecified initial multiplier
values. The nonlinearity makes it hard to construct an accurate
estimate of sensitivity. In addition this problem, the symmetry and
the underactuation induce numerical ill-conditioning. Therefore, in
order to apply the shooting method, it is important to compute the
sensitivities accurately, and the effects of the symmetry should be
taken into account.

In this paper, the attitude dynamics of a rigid body is described by
the structure-preserving Lie group variational integrator, and the
sensitivity is expressed in terms of a Lie algebra element. This
approach completely avoids any singularity in the attitude
representation, and the discrete flow respects the geometric
features. The resulting sensitivity derivatives are sufficiently
accurate for the shooting method. Furthermore, a simple numerical
approach is presented to eliminate the ill-conditioning caused by
the symmetry.

\paragraph*{Sensitivity derivatives} Taking a variation of the
discrete equations of motion and the multiplier equation using the
variational models, the linearized equations of motion and the
linearized multiplier equations can be written as
\begin{align*}
x_{k+1} & = A^{11}_k x_k + A^{12}_k \delta\lambda_k,\\
\delta\lambda_k & =A^{21}_{k+1} x_{k+1}+\parenth{A^{11}_{k+1}}^T
\delta\lambda_{k+1},
\end{align*}
where $x_k=[\zeta_k;\delta\Pi_k]\in\Re^6$, and matrices
$A^{ij}\in\Re^{6\times 6}$ are suitably defined. The solution of the
linear equations is given by
\begin{align*}
\begin{bmatrix}x_{N}\\\delta\lambda_{N}\end{bmatrix}
=
\begin{bmatrix}\Psi^{11} & \Psi^{12}\\\Psi^{21} &
\Psi^{22}\end{bmatrix}
\begin{bmatrix}x_{0}\\\delta\lambda_{0}\end{bmatrix},
\end{align*}
where $\Psi^{ij}\in\Re^{6\times 6}$. For the given two point
boundary value problem, the initial attitude and the initial angular
momentum are fixed, and the terminal multiplier is free. Thus, we
have the following sensitivity equation for the terminal attitude
and the terminal angular momentum with respect to the initial
multiplier;
\begin{align}
x_N =\Psi^{12}\delta\lambda_0.\label{eqn:Psi12}
\end{align}

\paragraph*{Avoiding numerical ill-conditioning}
The symmetry yields a conserved quantity by Noether's theorem, and
it causes a fundamental singularity in the sensitivity derivatives
for the two point boundary value problem. At each iteration, we
require the inverse of the sensitivity derivative represented by the
matrix $\Psi^{12}$ to update the initial multiplier to satisfy the
terminal boundary condition. However, this sensitivity matrix has a
theoretical rank deficiency of one since the vertical component of
the inertial angular momentum is conserved regardless of the initial
multiplier variation. Therefore, equation \refeqn{Psi12} is
numerically ill-conditioned.

Here we presents a simple numerical scheme to avoid the numerical
ill-conditioning caused by the symmetry. We decompose the
sensitivity derivative into symmetric parts and asymmetric parts.
Equation \refeqn{Psi12} is rewritten as
\begin{align}
\begin{bmatrix} \zeta_N \\ \delta\Pi_N \end{bmatrix}
=\begin{bmatrix} \Psi_1 & \Psi_2 \\ \Psi_3 &\Psi_4\end{bmatrix}
\begin{bmatrix} \delta\lambda^1_0 \\ \delta\lambda^2_0
\end{bmatrix},\label{eqn:Psi}
\end{align}
where $\Psi_i\in\Re^{3\times 3}$ are submatrices of $\Psi^{12}$.
Using the above equation and \refeqn{delRk}, the infinitesimal
variation of the inertial angular momentum is given by
\begin{align*}
\delta\pi_N & = \delta(R_N\Pi_N)=\delta R_N \Pi_N +
R_N\delta\Pi_N,\\
&=-R_NS(\Pi_N)\zeta_N + R_N\delta\Pi_N,\\
& =
-R_NS(\Pi_N)(\Psi_1\lambda^1_0+\Psi_2\lambda^2_0)+R_N(\Psi_3\lambda^1_0+\Psi_4\lambda^2_0).
\end{align*}
Now, the sensitivity derivative equation \refeqn{Psi} can be
rewritten in terms of the inertial angular momentum variation as
\begin{align}
\begin{bmatrix} \zeta_N \\ \delta\pi_N \end{bmatrix}
&=\begin{bmatrix} \Psi_1 & \Psi_2 \\ R_N(\Psi_3-S(\Pi_N)\Psi_1)
&R_N(\Psi_4-S(\Pi_N)\Psi_2)\end{bmatrix}
\begin{bmatrix} \delta\lambda^1_0 \\ \delta\lambda^2_0\end{bmatrix}.\label{eqn:Psi'}
\end{align}
\clearpage\newpage\noindent From the symmetry, the third component
of the inertial angular momentum variation is zero;
$\delta(\pi_N)_3=0$. Thus, the sixth row of the above matrix is
zero. (Numerical simulation in the later section shows that the norm
of the last row of the transformed sensitivity matrix is at the
level of $10^{-15}$.) Now, we find an update of the initial
multiplier by the pseudo-inverse of the $5\times 6$ matrix;
\begin{align}\label{eqn:Xi}
\delta\lambda_0 = \Xi^\dag x'_N=\Xi^T(\Xi\Xi^T)^{-1}x'_N,
\end{align}
where $\Xi\in\Re^{5\times 6}$ is composed of the first five rows of
the transformed sensitivity derivative in \refeqn{Psi'}, and
$x'_N=[\zeta_N;\delta(\pi_N)_1;\delta(\pi_N)_2]\in\Re^5$. This
approach removes the singularity in the sensitivity derivatives
completely, and the resulting optimal control problem is no longer
ill-conditioned. Numerical simulations show that the numerical
optimization procedure fails without this modification.

\paragraph*{Newton iteration} Using the decomposed sensitivity, an initial guess of the unspecified
initial conditions is iterated to satisfy the specified terminal
boundary conditions in the limit. Any type of Newton iteration can
be applied. We use a line search with backtracking algorithm,
referred to as Newton-Armijo iteration~\cite{bk:kelley}. The
procedure is summarized as follows.

\vspace*{0.1cm} {
\renewcommand{\theenumi}{\arabic{enumi}}
\renewcommand{\labelenumi}{\theenumi:}
\newcommand{\tab}{\hspace*{0.6cm}}
\hrule\vspace*{0.08cm}
\begin{enumerate}
\item Guess an initial multiplier $\lambda_0$.
\item Find $\Pi_k,R_k,\lambda_k^{1},\lambda_k^{2}$ using \refeqn{findf1}--\refeqn{updatelam}.
\item Compute the terminal B.C. error; $\mathrm{Error}=\norm{x'_N}$.
\item Set $\mathrm{Error}^t=\mathrm{Error},\;\; i=1$.
\item \textbf{while} $\mathrm{Error} > \epsilon_S$.
\item \tab Find a line search direction; $D=\Xi^\dag$.
\item \tab Set $c=1$.
\item \tab\textbf{while} $\mathrm{Error}^t > (1-2\alpha c)\mathrm{Error}$
\item \tab\tab Choose a trial multiplier $\lambda_0^{t}=\lambda_0^{}+c D
z_N$.
\item \tab\tab Find $\Pi_k^{},R_k^{},\lambda_k^{1},\lambda_k^{2}$ using \refeqn{findf1}--\refeqn{updatelam}.
\item \tab\tab Compute the error; $\mathrm{Error}^t=\norm{z_N^t}$.
\item \tab\tab Set $c=c/10,\;\; i=i+1$.
\item \tab\textbf{end while}
\item \tab Set $\lambda_0=\lambda_0^t$, $\mathrm{Error}=\mathrm{Error}^t$. (accept the trial)
\item \textbf{end while}
\end{enumerate}
\vspace*{0.08cm} \hrule} \vspace*{0.1cm} \noindent Here $i$ is the
number of iterations, and $\epsilon_S,\alpha\in\Re$ are stopping
criterion and a scaling factor, respectively. The outer loop finds a
search direction by computing the sensitivity derivatives, and the
inner loop performs a line search to find the largest step size
$c\in\Re$ along the search direction. The error in satisfaction of
the terminal boundary condition is determined at each inner
iteration.

\section{Numerical Examples}\label{sec:ne}

\begin{figure}[b]
    \centerline{\subfigure[Body (A)]{
    \includegraphics[width=0.45\columnwidth]{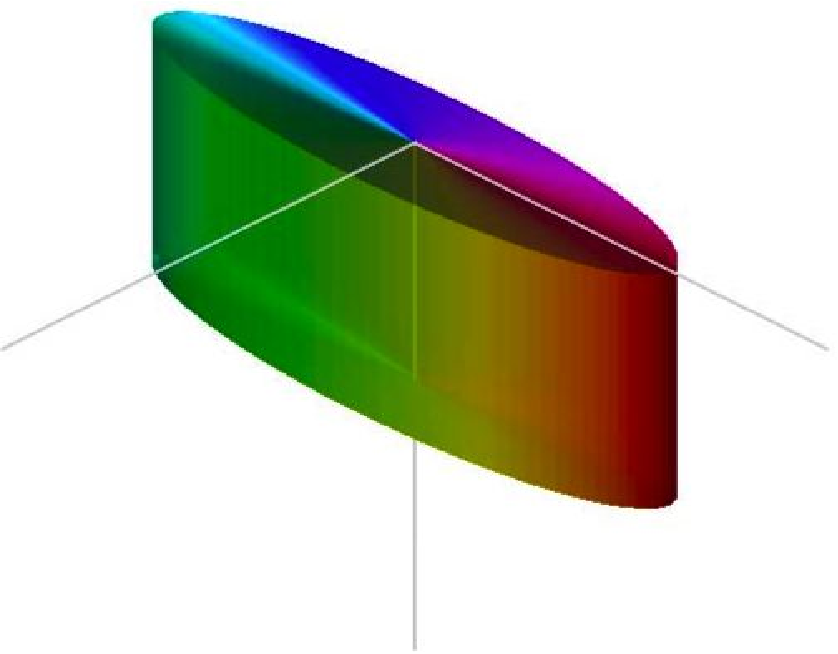}}
    \hspace*{1cm}
    \subfigure[Body (B)]{
    \includegraphics[width=0.45\columnwidth]{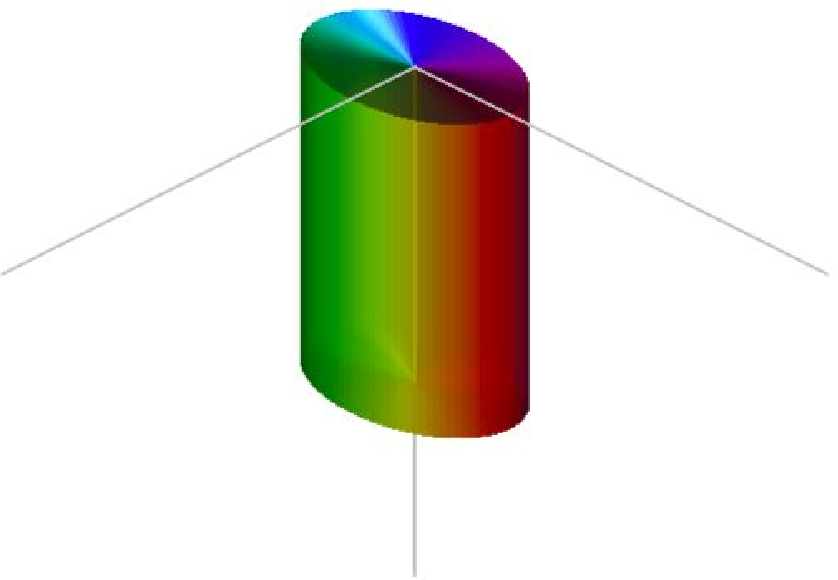}}
    }
    \caption{Elliptical cylinder}\label{fig:body}
\end{figure}

Numerical optimization results for the 3D pendulum are given. Two
elliptical cylinders, shown in \reffig{body}, are used as rigid
pendulum models. The properties are chosen as
\renewcommand{\labelenumi}{}
\begin{gather*}
\text{Body (A): } m=1,\,J=\mathrm{diag}[0.13,0.28,0.17],\rho=0.3
e_3.\\
\text{Body (B): } m=1,\,J=\mathrm{diag}[0.22,0.23,0.03],\rho=0.4
e_3.
\end{gather*}

Four cases are considered. Each maneuver is from a hanging
equilibrium to another hanging equilibrium with a rotation about the
vertical axis. The rotation angles are chosen as $90^\circ$ and
$180^\circ$. Since the vertical component of the angular momentum is
set to zero, the rotation is purely caused by the geometric phase
effect given in \refeqn{geo}. These problems are challenging in the
sense that the desired maneuvers are rotations about the gravity
direction, but the control input cannot directly generate any moment
about the gravity direction.

The corresponding boundary conditions are as follows.
\renewcommand{\theenumi}{\roman{enumi}}
\renewcommand{\labelenumi}{(\theenumi)}
\begin{enumerate}
\item Body (A), hanging equilibrium to hanging equilibrium with
$90^\circ$ yaw
\begin{gather*}
R_{0}=I_{3\times 3},\quad R_{N}^d=\begin{bmatrix}0&-1&0\\1&0&0\\0&0&1\end{bmatrix},\\
\Pi_0=0_{3\times 1},\quad\Pi_N^d=0_{3\times 1}.
\end{gather*}
\item Body (A), hanging equilibrium to hanging equilibrium with
$180^\circ$ yaw
\begin{gather*}
R_{0}=I_{3\times 3},\quad R_{N}^d=\mathrm{diag}[-1,-1,1],\\
\Pi_0=0_{3\times 1},\quad\Pi_N^d=0_{3\times 1}.
\end{gather*}
\item Body (B), hanging equilibrium to hanging equilibrium with
$90^\circ$ yaw
\begin{gather*}
R_{0}=I_{3\times 3},\quad R_{N}^d=\begin{bmatrix}0&-1&0\\1&0&0\\0&0&1\end{bmatrix},\\
\Pi_0=0_{3\times 1},\quad\Pi_N^d=0_{3\times 1}.
\end{gather*}
\item Body (B), hanging equilibrium to hanging equilibrium with
$180^\circ$ yaw
\begin{gather*}
R_{0}=I_{3\times 3},\quad R_{N}^d=\mathrm{diag}[-1,-1,1],\\
\Pi_0=0_{3\times 1},\quad\Pi_N^d=0_{3\times 1}.
\end{gather*}
\end{enumerate}

\begin{table}\normalsize\selectfont
\begin{center}
\begin{threeparttable}
\caption{Optimization results}\label{tab:opt}
\renewcommand{\arraystretch}{1.1}
\begin{tabular}{c|>{$}c<{$}>{$}c<{$}>{$}c<{$}|c}\hline\hline
Case  & \mathcal{J} & $\norm{\mathrm{logm}(R_N^{d,T}R_N)}$&
$\norm{\Pi_N^d-\Pi_N}$ & $\Delta T$\\\hline
(i)   & 5.91 & 2.30\times 10^{-14} & 1.34\times 10^{-14} & 2.72\\
(ii)  & 7.32 & 4.80\times 10^{-15} & 1.66\times 10^{-14} & 5.25\\
(iii) & 1.73 & 1.22\times 10^{-15} & 6.55\times 10^{-14} & 4.09\\
(iv)  & 3.37 & 3.06\times 10^{-14} & 3.04\times 10^{-14} &
5.05\\\hline\hline
\end{tabular}
\begin{tablenotes}
\item $\Delta T$: Simulation running time in Intel Pentium M 740 1.73GHz
processor (min.)
\end{tablenotes}
\end{threeparttable}
\end{center}
\end{table}

The optimal control results are given in Table \ref{tab:opt}, where
the optimized performance index, the error in satisfaction of the
terminal boundary condition, and the simulation running time are
shown for each case. The terminal error is at the level of machine
precision, and the simulation time is about 5 minutes.

\begin{figure}
    \centerline{\subfigure[Attitude Maneuver]{
    \includegraphics[width=0.5\columnwidth]{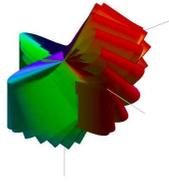}
    }
    \hfill
    \subfigure[Geometric Phase]{
    \includegraphics[width=0.48\columnwidth]{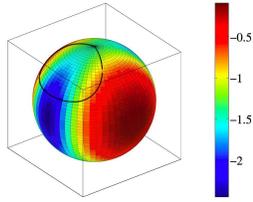}}
    }
    \centerline{\subfigure[Control Input]{
    \includegraphics[width=0.45\columnwidth]{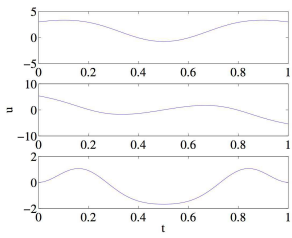}}
    \hfill
    \subfigure[Convergence Rate]{
    \includegraphics[width=0.47\columnwidth]{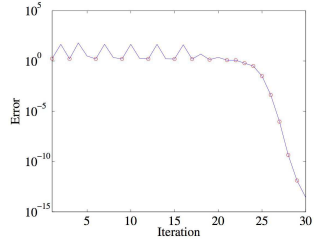}}
    }
    \caption{(i) Body A, hanging equilibrium to hanging equilibrium with
$90^\circ$ yaw}\label{fig:i}
\end{figure}

\begin{figure}
    \centerline{\subfigure[Attitude Maneuver]{
    \includegraphics[width=0.5\columnwidth]{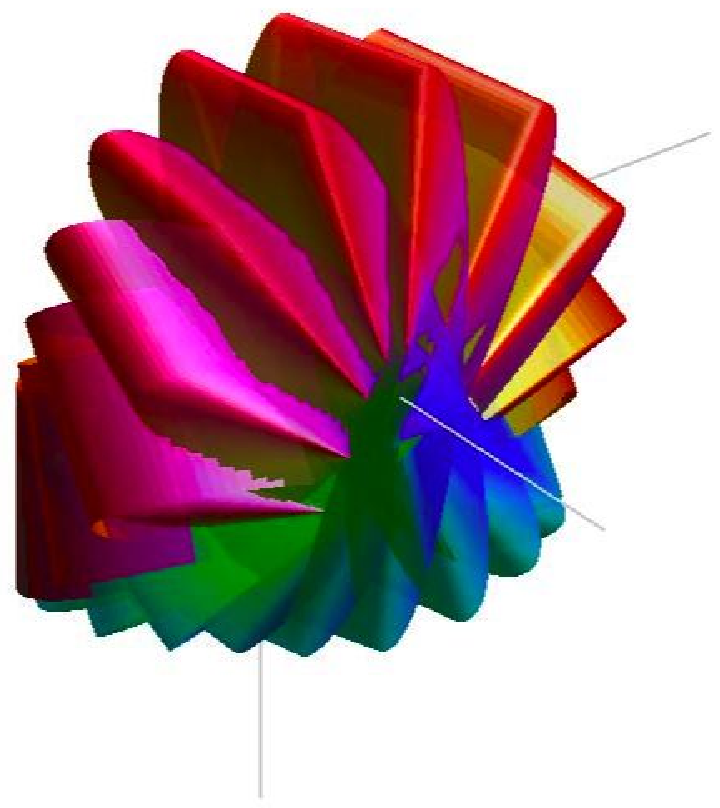}
    }
    \hfill
    \subfigure[Geometric Phase]{
    \includegraphics[width=0.48\columnwidth]{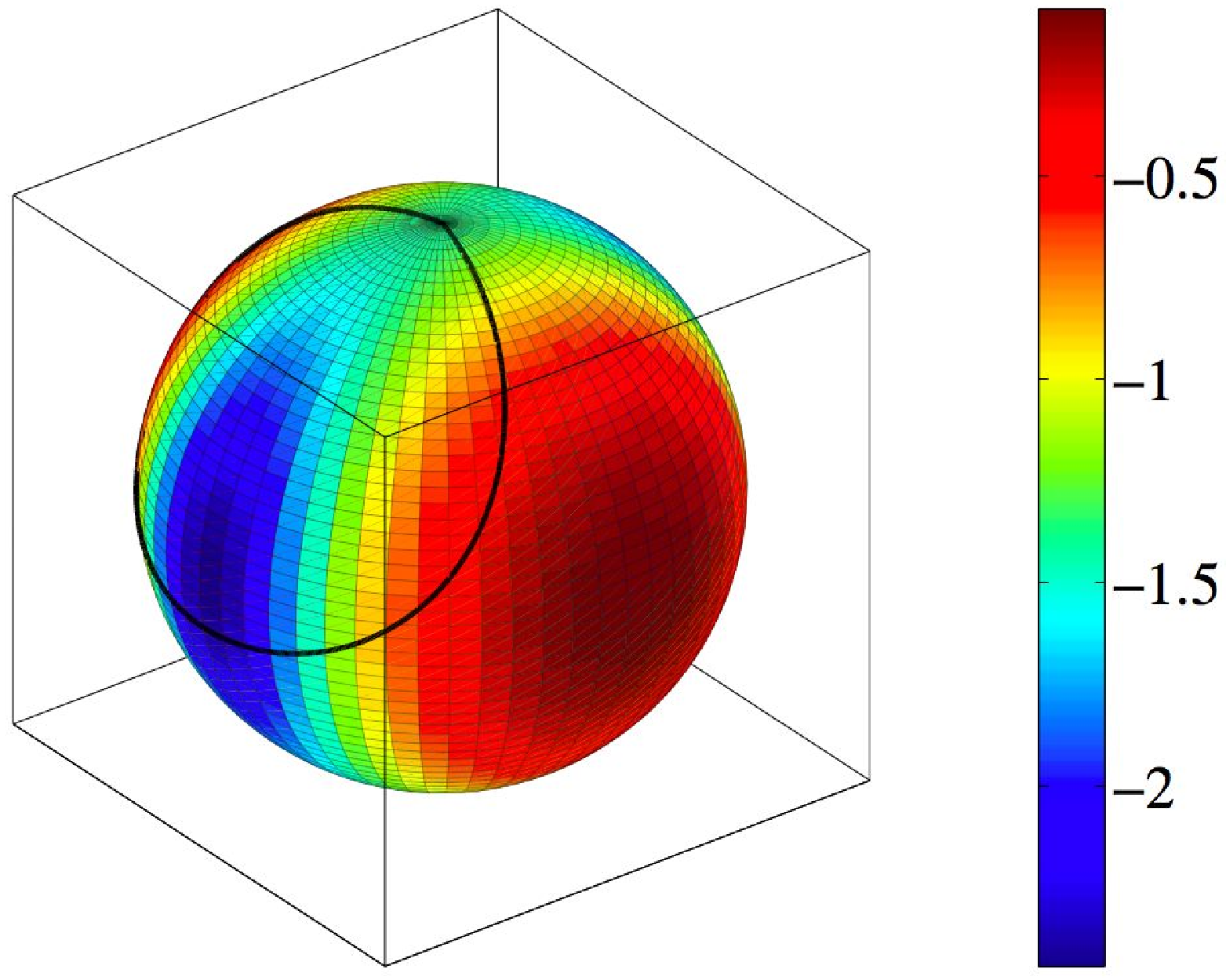}}
    }
    \centerline{\subfigure[Control Input]{
    \includegraphics[width=0.45\columnwidth]{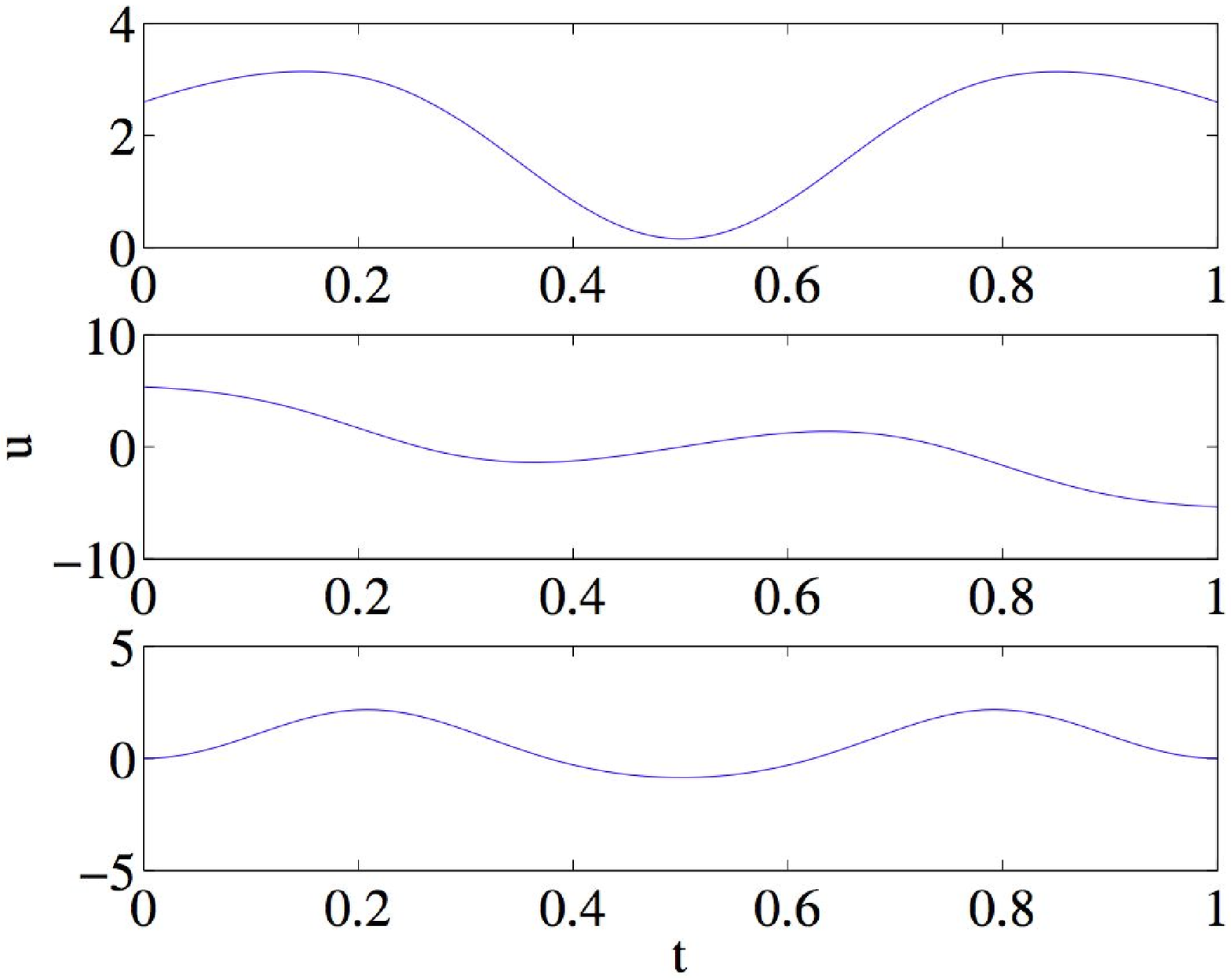}}
    \hfill
    \subfigure[Convergence Rate]{
    \includegraphics[width=0.47\columnwidth]{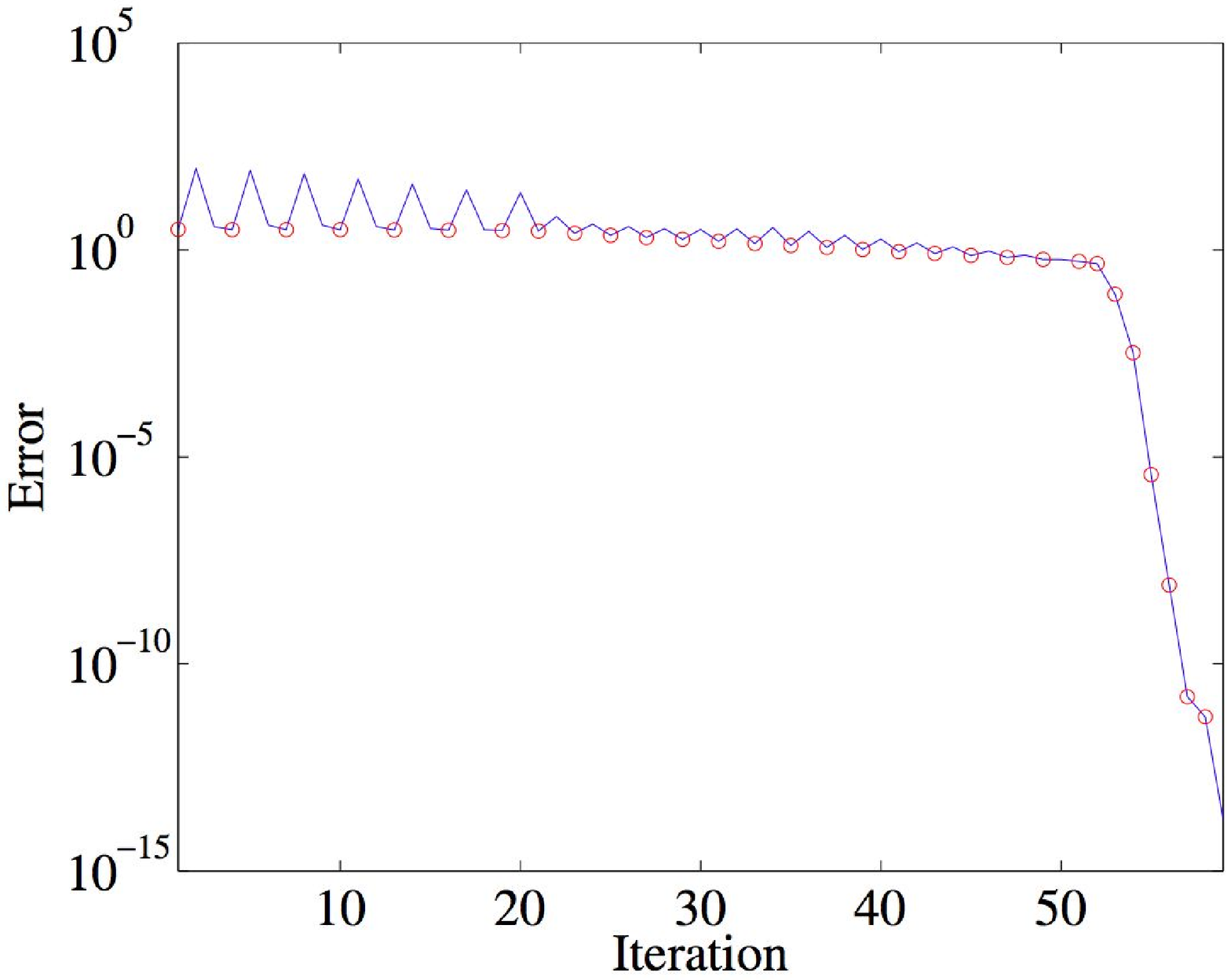}}
    }
    \caption{(ii). Body A, hanging equilibrium to hanging equilibrium with
$180^\circ$ yaw}\label{fig:ii}
\end{figure}

\begin{figure}
    \centerline{\subfigure[Attitude Maneuver]{
    \includegraphics[width=0.52\columnwidth]{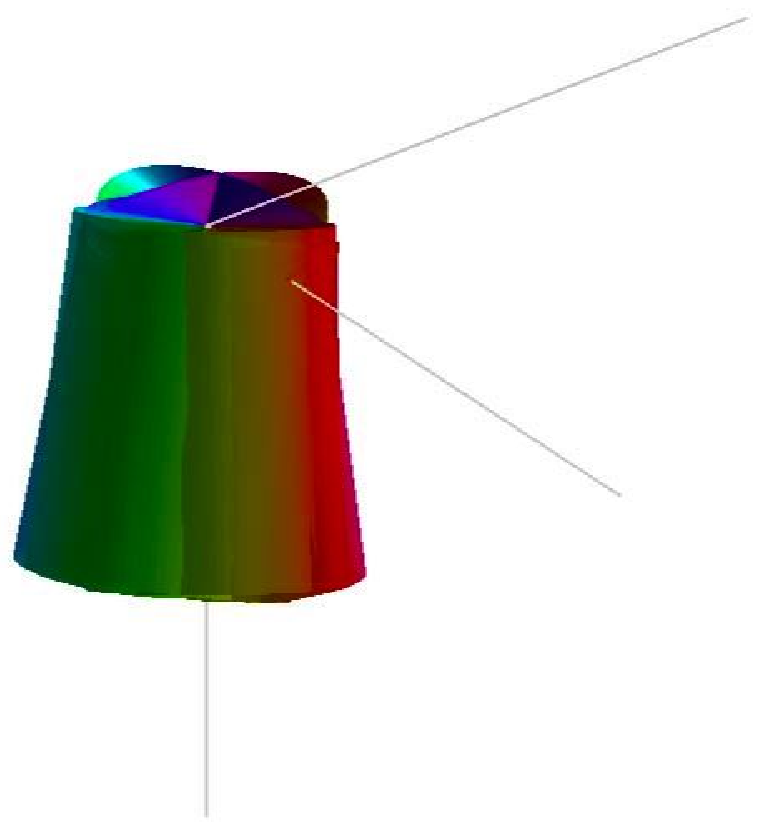}
    }
    \hfill
    \subfigure[Geometric Phase]{
    \includegraphics[width=0.50\columnwidth]{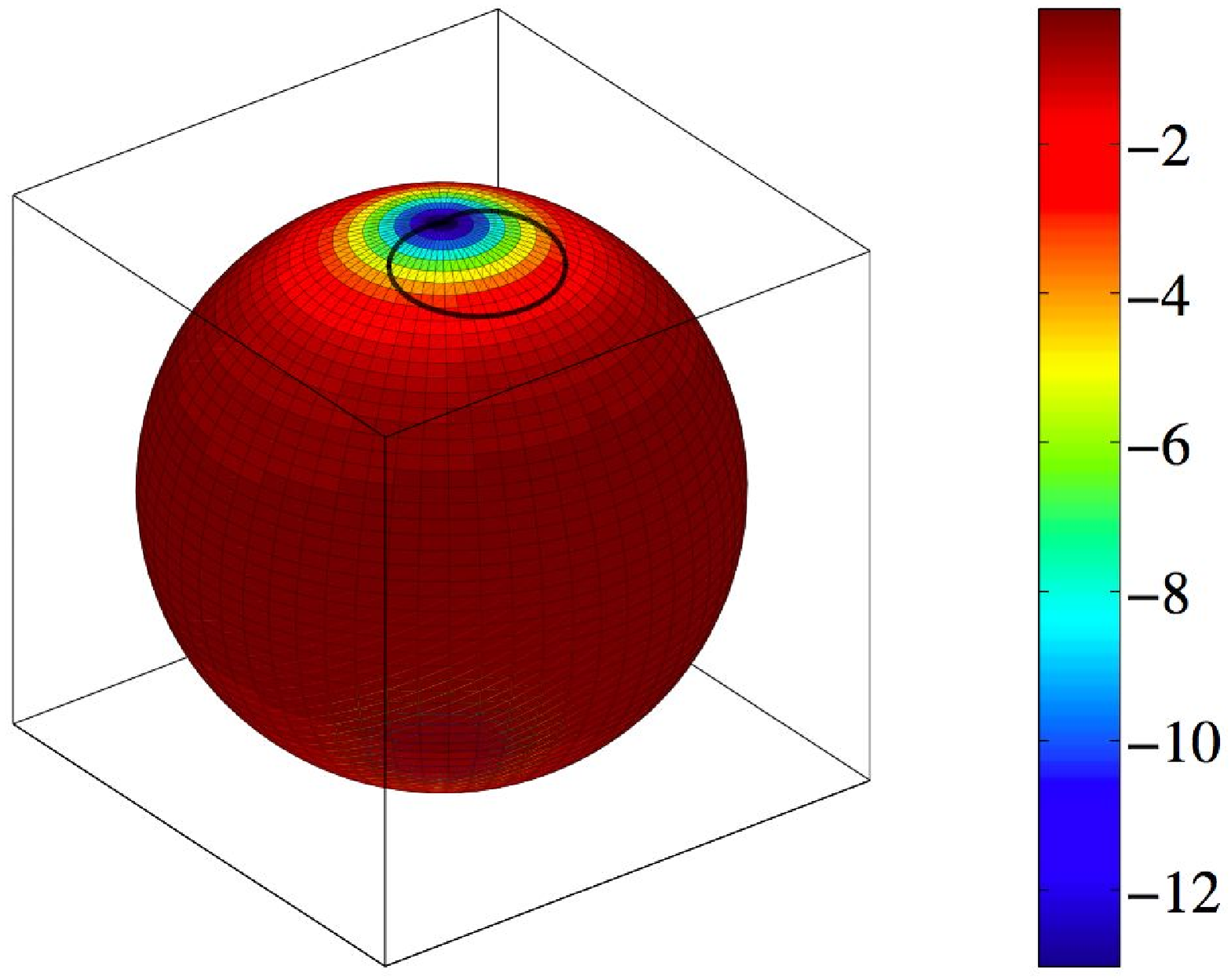}}
    }
    \centerline{\subfigure[Control Input]{
    \includegraphics[width=0.45\columnwidth]{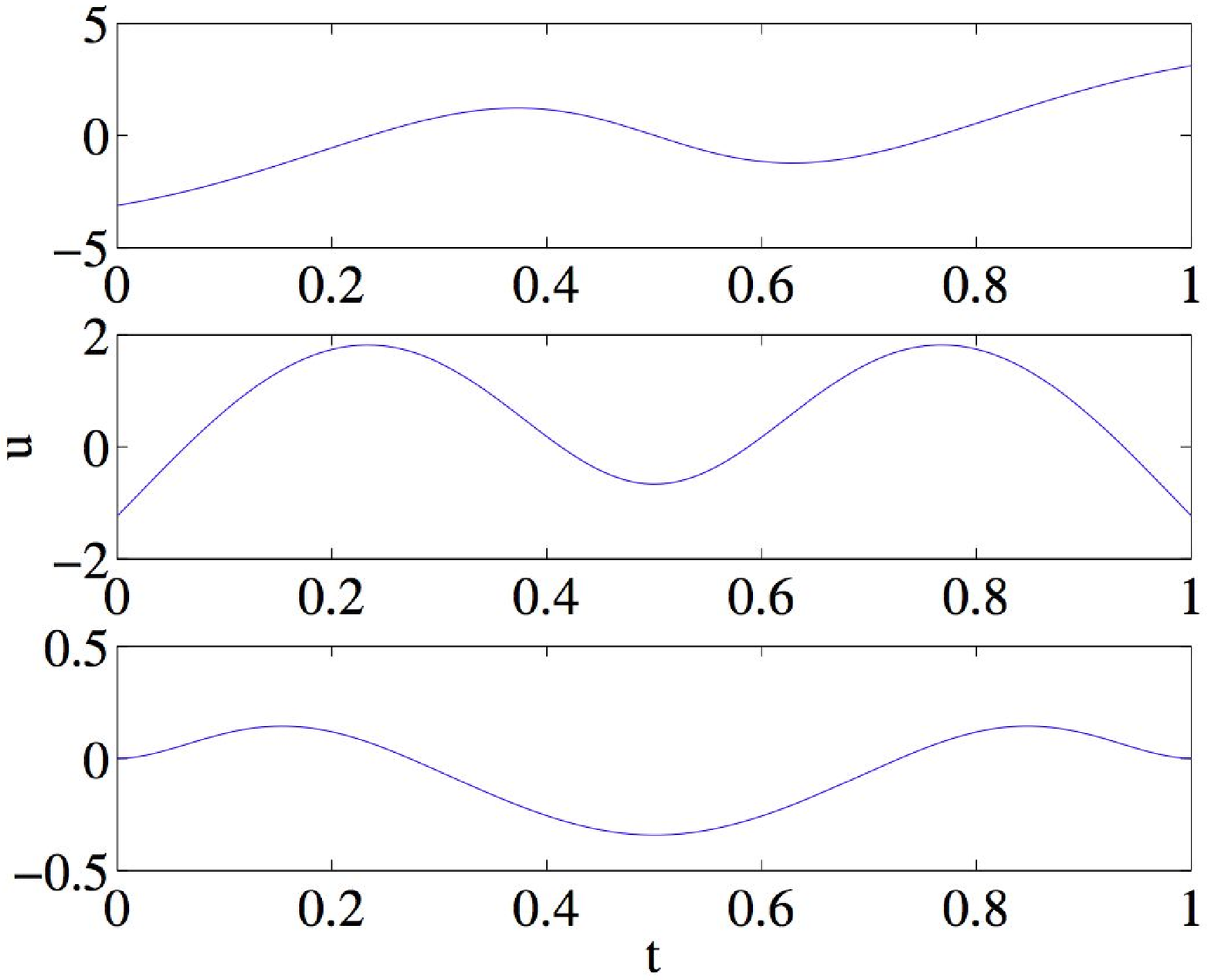}}
    \hfill
    \subfigure[Convergence Rate]{
    \includegraphics[width=0.47\columnwidth]{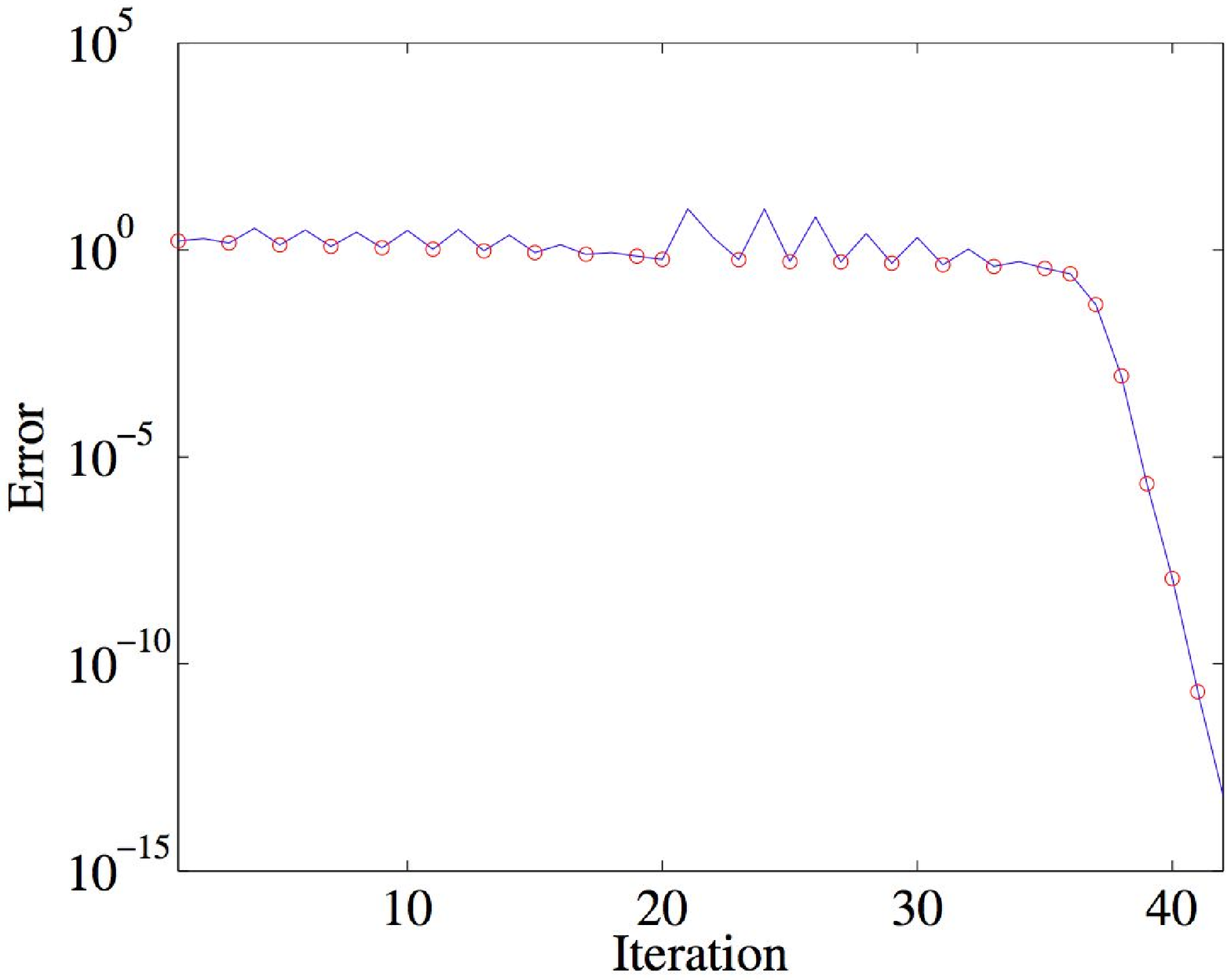}}
    }
    \caption{(iii). Body B, hanging equilibrium to hanging equilibrium with
$90^\circ$ yaw}\label{fig:iii}
\end{figure}

\begin{figure}
    \centerline{\subfigure[Attitude Maneuver]{
    \includegraphics[width=0.51\columnwidth]{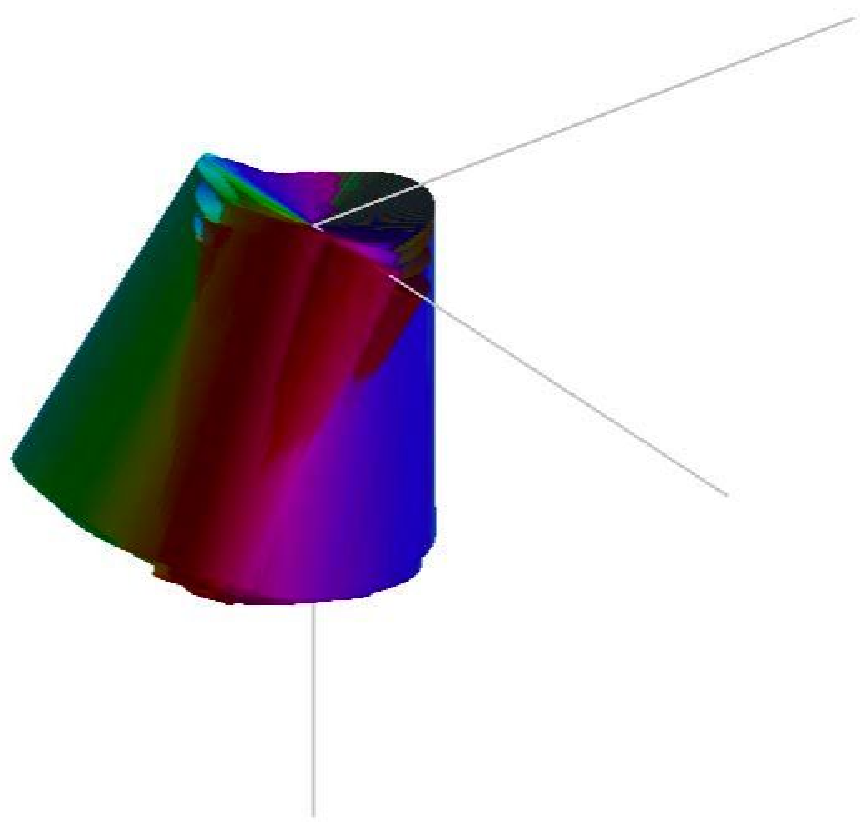}
    }
    \hfill
    \subfigure[Geometric Phase]{
    \includegraphics[width=0.5\columnwidth]{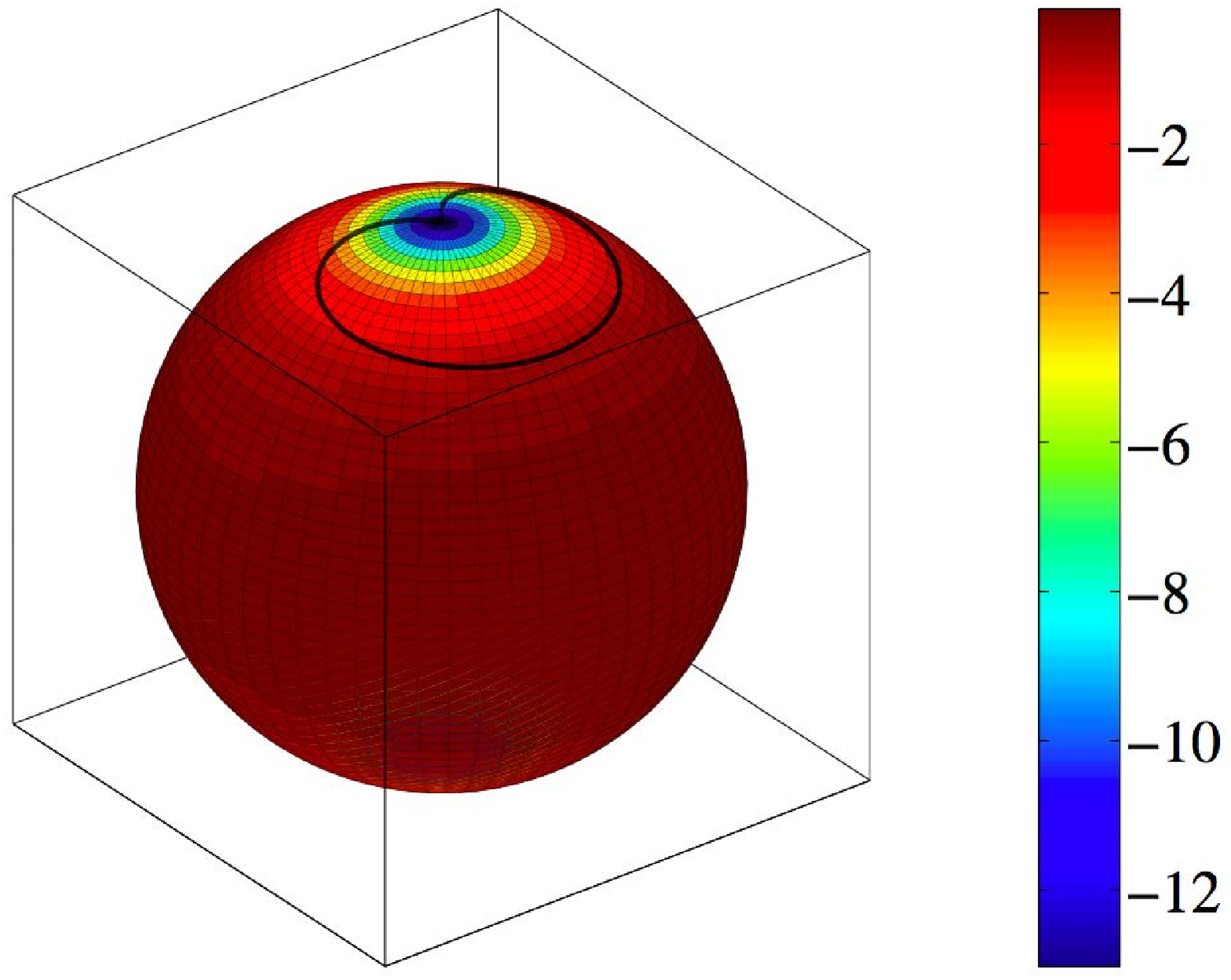}}
    }
    \centerline{\subfigure[Control Input]{
    \includegraphics[width=0.45\columnwidth]{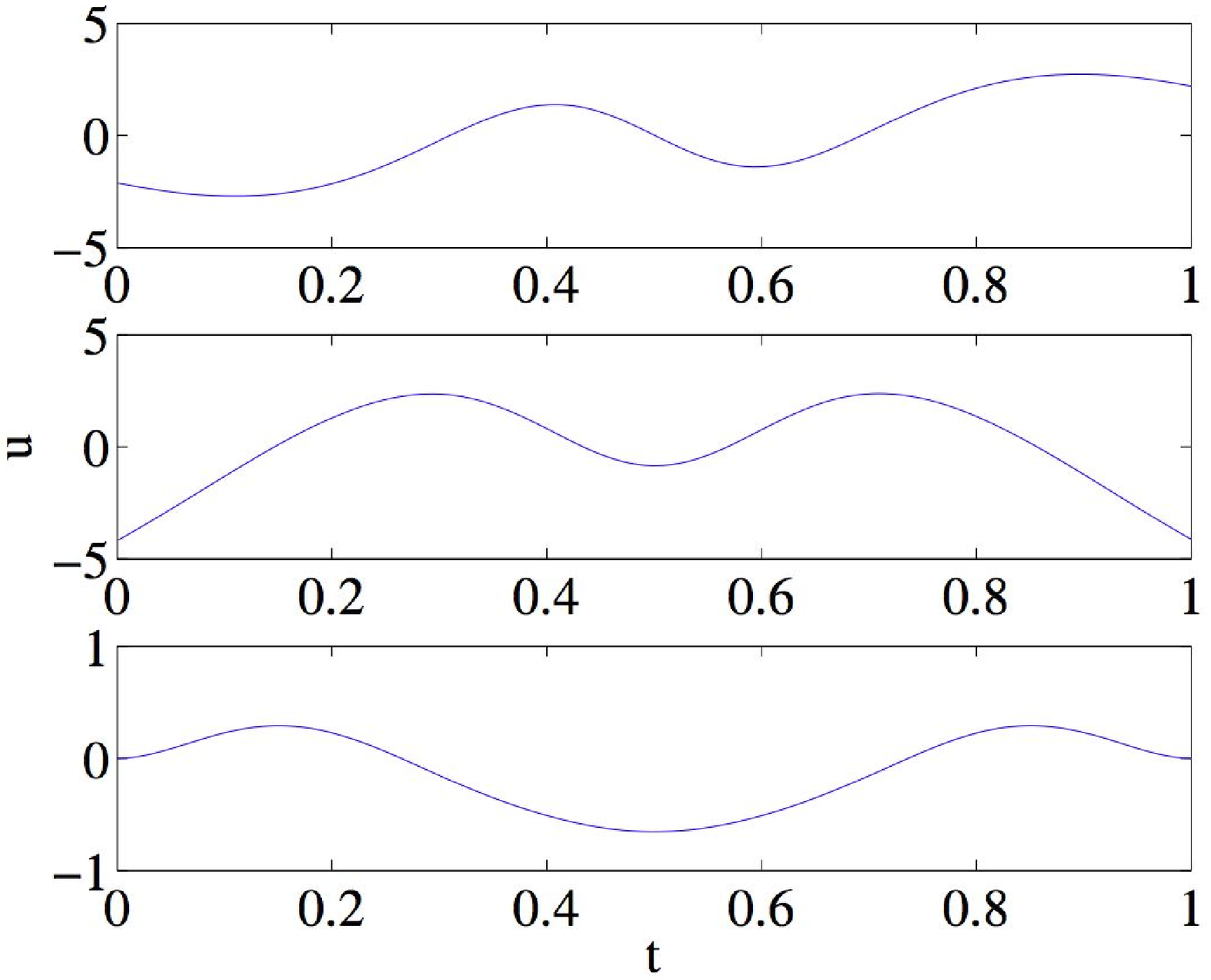}}
    \hfill
    \subfigure[Convergence Rate]{
    \includegraphics[width=0.47\columnwidth]{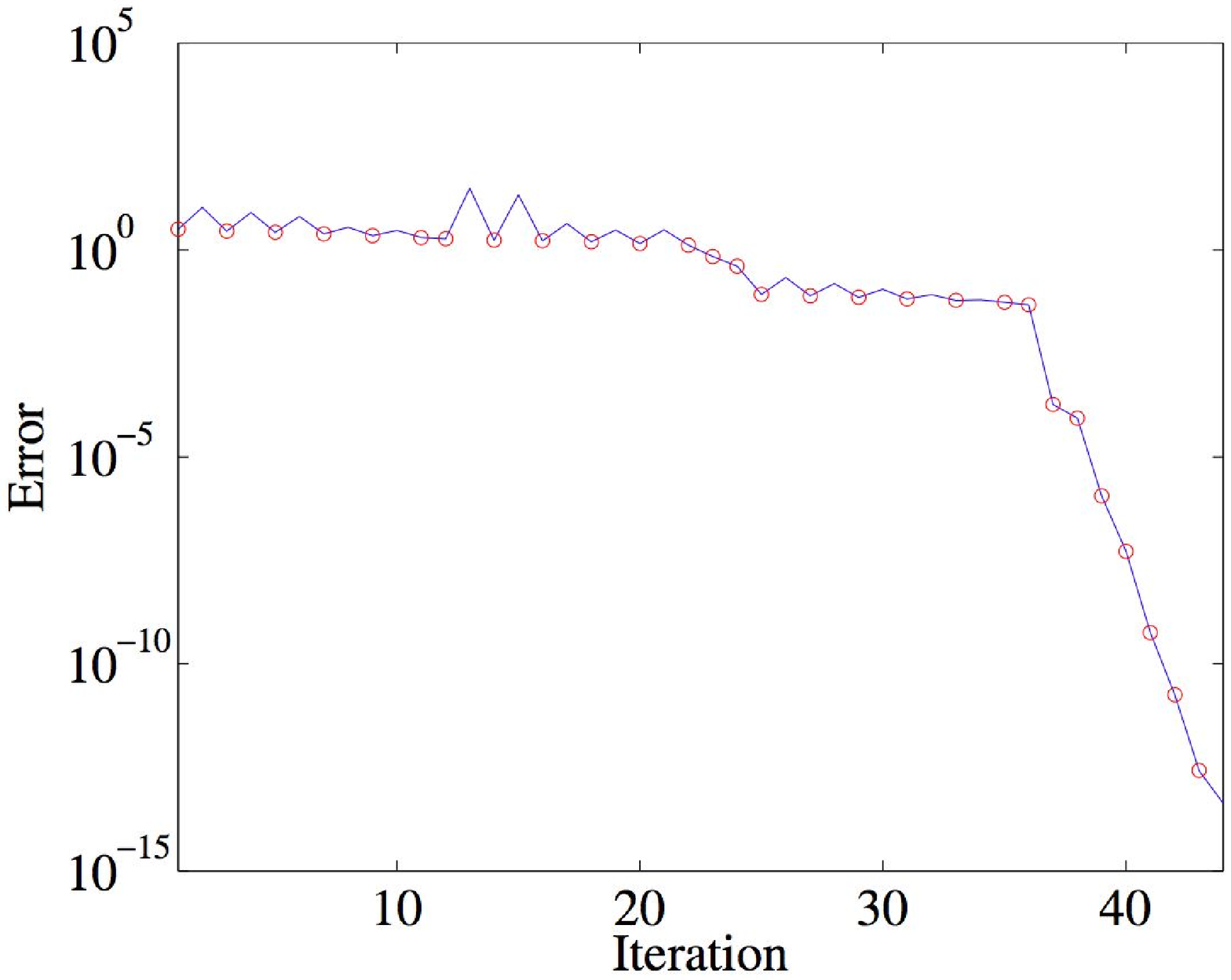}}
    }
    \caption{(iv). Body B, hanging equilibrium to hanging equilibrium with
$180^\circ$ yaw}\label{fig:iv}
\end{figure}

Figures \ref{fig:i}--\ref{fig:iv} show snapshots of the attitude
maneuvers, reduced trajectory of $\Gamma=R^T e_3$ on a sphere,
control input history, and convergence rate. (A simple animation for
the attitude maneuver can be seen at
\url{http://www.umich.edu/~tylee}.)

The convergence rate figures show violation of the terminal boundary
condition according to the number of iterations in a logarithm
scale. Red circles denote outer iterations in Newton-Armijo
iteration to compute the sensitivity derivatives. For all cases, the
initial guesses of the unspecified initial multiplier are
arbitrarily chosen. The error in satisfaction of the terminal
boundary condition converges quickly to machine precision after the
solution is close to the local minimum at around 50th iteration.
These convergence results are consistent with the quadratic
convergence rates expected of Newton methods with accurately
computed gradients. The condition number of the decomposed
sensitivity derivative given at \refeqn{Xi} varies from $10^0$ to
$10^5$. If the sensitivity derivative is not decomposed, then the
condition numbers are at the level of $10^{19}$, and the numerical
iterations fail.

The numerical examples presented in this paper show excellent
numerical convergence properties. This is because the proposed
computational algorithms on \SO are geometrically exact and
numerically accurate. In addition, the algorithm incorporates a
modification that eliminates the singularity caused by the symmetry.

We interpret the optimization results using the geometric phase
formula given by \refeqn{geo}. For given initial conditions, the
vertical component of the initial angular momentum is zero. Thus,
the rotation about the vertical axis is purely caused by the
geometric phase. Since the geometric phase is determined by a
surface integral on $\S^2$ whose boundary is the reduced trajectory
$\Gamma$, it is more efficient for the reduced trajectory to enclose
the area at which the absolute value of the integrand of
\refeqn{geo} is maximized.

In each subfigure (b) of Figures \ref{fig:i}--\ref{fig:iv}, the
infinitesimal geometric phase per unit area is shown by color
shading. The reduced trajectory, which represents the gravity
direction in the body fixed frame, is shown by a solid line. The
north pole of the sphere corresponds to the hanging equilibrium
manifold, and the reduced trajectory starts and ends at the same
north pole for the given boundary conditions.

Comparing Figures 2(b), 3(b) with Figure 4(b), 5(b), it can be seen
that Body (A) and Body (B) have different geometric phase
characteristics. This is caused by the fact that the geometric phase
depends on the moment of inertia of the body. For Body (A), the
absolute value of the infinitesimal geometric phase is maximized at
a point on the equator, and for Body (B), it is maximized at the
north pole. We see that the optimized reduced trajectories try to
enclose those points.

As a result, the optimized attitude maneuver of Body (A) is
distinguished from that of Body (B). The attitude maneuver of Body
(A) is relatively more aggressive than that of Body (B) since the
reduced trajectory passes near the equator corresponding to a
horizontal position. Body (B) does not have to move far away from
the hanging equilibrium since the infinitesimal geometric phase is
maximized at that point. The resulting attitude maneuver is
relatively benign.

\bibliography{opt}
\bibliographystyle{IEEEtran}

\clearpage\newpage

\end{document}